\magnification=\magstep1 
\baselineskip=14pt 
\def\mapright#1{\smash{\mathop{\longrightarrow}\limits^{#1}}} 
\def\frac#1#2{{#1\over #2}} 

% Reference macros. Sets up a count as a counter for refs, resets it to 0.
\newcount\refno\refno=0
% References are typed by putting \ref at the start of the line, the text
% of the reference, and a blank line afterwards.  After the blank line
% following a reference to (say) Ladner, put
%     \newcount\ladner\ladner=\refno
% Then in your document, type [\the\ladner] to obtain [XX]
% At the end of the list of references, the line :
%            \immediate\closeout\reffile
% should appear to close the references file, and then
%            \input TempReferences
% will read them in later.
\chardef\other=12
\newwrite\reffile
\immediate\openout\reffile=TempReferences
\outer\def\ref{\par\medbreak\global\advance\refno by 1
  \immediate\write\reffile{}
  \immediate\write\reffile{\noexpand\item{[\the\refno]}}
  \copytoblankline}
\def\copytoblankline{\begingroup\setupcopy\copyref}
\def\setupcopy{\def\do##1{\catcode`##1=\other}\dospecials
  \catcode`\|=\other \obeylines}
{\obeylines \gdef\copyref#1
  {\def\next{#1}%
  \ifx\next\empty\let\next=\endgroup %
  \else\immediate\write\reffile{\next} \let\next=\copyref\fi\next}}

\ref ALLEN, S.D.; SINCLAIR, A.M.; SMITH, R.R., The ideal structure of the Haagerup tensor product of C$^*$-algebras, {\sl J. reine angew. Math.}, 442 (1993) 111-148.

\newcount\ASS\ASS=\refno

\ref ARCHBOLD, R.J., Topologies for primal ideals, {\sl J. London
Math. Soc.}, (2) 36 (1987) 524-542.

\newcount\RJA\RJA=\refno

\ref ARCHBOLD, R.J.; KANIUTH, E.; SCHLICHTING, G.; SOMERSET, D.W.B., Ideal spaces of the Haagerup tensor product of C$^*$-algebras, {\sl Internat. J. Math.}, 8 (1997) 1-29.

\newcount\AKSS\AKSS=\refno

\ref BAGGETT, L.W., A separable group having a discrete dual space is compact, {\sl J. Funct. Anal}, 10 (1972) 131-148.

\newcount\Ba\Ba=\refno

\ref BARNES, B.A., Ideal and representation theory of the $L^1$-algebra of a group with polynomial growth, {\sl Colloq. Math.}, 45 (1981) 301-315.

\newcount\Bar\Bar=\refno

\ref BECKHOFF, F., Topologies on the space of ideals of a Banach algebra, {\sl Stud. Math.}, 115 (1995) 189-205.

\newcount\Be\Be=\refno

\ref BECKHOFF, F., Topologies of compact families on the ideal space of a Banach algebra, {\sl Stud. Math.}, 118 (1996) 63-75.

\newcount\Bec\Bec=\refno

\ref BECKHOFF, F.,  Topologies on the ideal space of a Banach algebra and spectral synthesis, {\sl Proc. Amer. Math. Soc.}, 125 (1997) 2859-2866.

\newcount\Beck\Beck=\refno

%\ref BLECHER, D.P.,  Geometry of the tensor product of C$^*$-algebras, {\sl Math. Proc. Camb. Phil. Soc.}, 104 (1988) 119-127.

%\newcount\Ble\Ble=\refno

%\ref BLECHER, D.P., Tensor products which do not preserve operator algebras, {\sl Math. Proc. Camb. Phil. Soc.}, 108 (1990) 395-403.

%\newcount\Blec\Blec=\refno

%\ref BOHNENBLUST, H.F.; KARLIN, S., Geometrical properties of the unit sphere of Banach algebras, {\sl Ann. Math.}, (2) 62 (1955) 217-229.

%\newcount\BK\BK=\refno

\ref BOIDOL, J.; LEPTIN, H.; SCH\"URMAN J.; VAHLE, D., R\"amme primitiver Ideale von Gruppenalgebren, {\sl Math. Ann.} 236 (1978) 1-13.

\newcount\BLSV\BLSV=\refno

\ref BONSALL, F.F.; DUNCAN, J., {\sl Complete Normed Algebras}, Springer-Verlag, New York, 1973.

\newcount\BD\BD=\refno

%\ref CHRISTENSEN, E.; SINCLAIR, A.M., A survey of completely bounded operators, {\sl Bull. London Math. Soc.}, 21 (1989) 417-448.

%\newcount\CS\CS=\refno

%\ref COHN, P.M., {\sl Algebra, Vol. 3} (2nd edn.), Wiley, Chichester, 1991.

%\newcount\Cohn\Cohn=\refno

%\ref DALES, H.G., 

%\newcount\Dal\Dal=\refno

%\ref DIXON, P.G., Simple radical Banach algebras and the Somerset condition.

%\newcount\Dix\Dix=\refno

\ref EYMARD, P., L'alg\`ebre de Fourier d'un groupe localement compact, {\sl Bull. Soc. Math. France}, 92 (1964) 181-236.

\newcount\Ey\Ey=\refno

\ref FEINSTEIN, J.F.; SOMERSET, D.W.B., A note on ideal spaces of Banach algebras, {\sl Bull. London Math. Soc.}, 30 (1998) 611-617.

\newcount\FS\FS=\refno

%\ref FEINSTEIN, J.F.; SOMERSET, D.W.B., Strong regularity for uniform algebras, 
%in {\sl Contemporary Math., 232}, `Proceedings of the 3rd Function Spaces Conference, Edwardsville, Illinois, May 1998', ed. K. Jarosz, pp. 139-149, Amer. Math. Soc., Rhode Island, 1999.

%\newcount\Fein\Fein=\refno

\ref FEINSTEIN, J.F.; SOMERSET, D.W.B., Spectral synthesis for Banach algebras, II, preprint (Nottingham University), 1999.

\newcount\fes\fes=\refno

\ref FORREST, B., Amenability and bounded approximate identities in ideals of $A(G)$, {\sl Illinois Math. J.}, 34 (1990) 1-25.

\newcount\For\For=\refno

%\ref GIERZ, G.; HOFMANN, K.H.; KEIMEL, K.; LAWSON, J.; MISLOVE, M.;
%SCOTT, D.S., {\sl A Compendium of Continuous Lattices}, Springer-Verlag,
%New York, 1980.

%\newcount\Comp\Comp=\refno

\ref GRAHAM, C.C.;\hfil McGEHEE, O.C.,\hfil{\sl Essays in Commutative 
Harmonic Analysis}, 
Springer-Verlag, New York, 1979.

\newcount\GM\GM=\refno

\ref GROSSER, S.; MOSKOWITZ, On central topological groups, {\sl Trans. Amer. Math. Soc.}, 127 (1967) 317-340.

\newcount\GM\GM=\refno

\ref HAUENSCHILD, W.; KANIUTH, E.; KUMAR, A., Ideal structure of Beurling algebras on [FC]$^-$ groups, {\sl J. Funct. Anal.}, 51 (1983) 213-228.

\newcount\HKK\HKK=\refno

\ref HAUENSCHILD, W.; LUDWIG, J., The injection and the projection theorem for spectral sets, {\sl Monatsh. Math.}, 92 (1981) 167-177.

\newcount\HaLu\HaLu=\refno

\ref HENRIKSEN, M.; KOPPERMAN, R.; MACK, J.; SOMERSET, D.W.B., Joincompact
spaces, continuous lattices, and C$^*$-algebras, {\sl Algebra Universalis}, 38 (1997) 289-323.

\newcount\HKMS\HKMS=\refno

%\ref HULANICKI, A.?
%
%\newcount\Hul\Hul=\refno

\ref KANIUTH, E., Primitive ideal spaces of groups with relatively compact conjugacy classes, {\sl Arch. Math.}, 32 (1979) 16-24.

\newcount\Kant\Kant=\refno

%\ref KANIUTH, E., $*$-regularity of locally compact groups, {\sl Probability Measures on Groups, VII, Oberwolfach (1983)}, 235-240, Lecture Notes in Math. 1064, Springer, Berlin, 1984.

%\newcount\Kaniu\Kaniu=\refno

\ref KANIUTH, E., The Helson-Reiter theorem for a class of nilpotent discrete groups, {\sl Math. Proc. Camb. Phil. Soc.}, 122 (1997) 95-103.

\newcount\Kan\Kan=\refno

\ref KANIUTH, E.; LAU, A.T., Spectral synthesis for $A(G)$ and subspaces of $VN(G)$, preprint (Paderborn), 1999.

\newcount\KL\KL=\refno

%\ref KANIUTH, E.; TAYLOR, K.F., Kazhdan constants and the dual space topology, {\sl Math. Ann.}, 293 (1992) 495-508.

%\newcount\KaTa\KaTa=\refno

\ref KAPLANSKY, I., The structure of certain operator algebras, {\sl Trans. Amer. Math. Soc.}, 70 (1951) 219-255.

\newcount\Kap\Kap=\refno

\ref KELLEY, J.L., {\sl General Topology}, Van Nostrand, Princeton, 1955.

\newcount\Kel\Kel=\refno

%\ref KOPPERMAN, R., Asymmetry and duality in topology, {\sl Gen. Top. and Appl.}, 66 (1995) 1-39.

%\newcount\Ko\Ko=\refno

%\ref LEPTIN, H.; \"Uber die schwache Wienersche Eigenschaft f\"ur Liesche Gruppen, {\sl Archiv. Math.}, 35 (1980) 144-148.

%\newcount\Lep\Lep=\refno

\ref LEPTIN, H.; Structure of $L^1(G)$ for locally compact groups, in {\sl Operator Algebras and Group Representations}, Vol. II, (Proceedings of Neprun Conference 1980), Pitman, London, 1984.

\newcount\Lept\Lept=\refno

\ref LEPTIN, H.; POGUNTKE, D., Symmetry and nonsymmetry for locally compact groups, {\sl J. Funct. Anal.}, 33 (1979) 119-134.

\newcount\LeP\LeP=\refno

\ref LUDWIG, J., A class of symmetric and a class of Wiener group algebras, {\sl J. Funct. Anal.}, 31 (1979) 187-194.

\newcount\Lu\Lu=\refno

\ref LUDWIG, J., Polynomial growth and ideals in group algebras, {\sl Manuscripta. Math.} 30 (1980) 215-221.

\newcount\Lud\Lud=\refno

\ref PALMER, T.W., Classes of nonabelian noncompact, locally compact groups, {\sl Rocky Mountain J. Math}, 8 (1978) 683-741.

\newcount\Palm\Palm=\refno

%\ref PALMER, T.W., {\sl Banach Algebras and the General Theory of $^*$-Algebras}, Vol. 1, C.U.P., New York, 1994.

%\newcount\Pal\Pal=\refno

%\ref PAULSEN, V.I.; SMITH, R.R., Multilinear maps and tensor norms on operator systems, {\sl J. Funct. Anal.} 73 (1987) 258-276.

%\newcount\PS\PS=\refno

%\ref RICKART, C.E., On spectral permanence for certain Banach algebras, {\sl Proc. Amer. Math. Soc.}, 4 (1953) 191-196.

%\newcount\Rick\Rick=\refno

%\ref RICKART, C.E., {\sl General Theory of Banach Algebras}, Van Nostrand, London, 1960.

%\newcount\Ri\Ri=\refno

\ref ROBERTSON, L., A note on the structure of Moore groups, {\sl Bull. Amer. Math. Soc.}, 75 (1969) 594-598.

\newcount\Robe\Robe=\refno

\ref RUDIN, W, {\sl Fourier Analysis on Groups}, Interscience, New York, 1962.

\newcount\Rud\Rud=\refno

\ref RUNDE, V., Intertwining operators over $L_1(G)$ for $G\in [{\rm PG}]\cap
[{\rm SIN}]$, {\sl Math. Z.} 221 (1996) 495-506.

\newcount\Run\Run=\refno

\ref SOMERSET, D.W.B., Spectral synthesis for Banach algebras, {\sl Quart. J. Math. Oxford}, (2) 49 (1998) 501-521.

\newcount\Syn\Syn=\refno

\ref SOMERSET, D.W.B., Ideal spaces of Banach algebras, {\sl Proc. London Math. Soc.}, (3) 78 (1999) 369-400.

\newcount\Id\Id=\refno

%\ref TAKESAKI, M., {\sl Theory of Operator Algebras, I}, Springer-Verlag, New York, 1979.

%\newcount\Tak\Tak=\refno

%\ref ZELAZKO, W., {\sl Banach Algebras}, Elsevier, Warsaw, 1973.

%\newcount\Zel\Zel=\refno

\ref ZELMANOV, E.I., On periodic compact groups, {\sl Israel J. Math}, 77 (1992) 83-95.

\newcount\Zelm\Zelm=\refno

\immediate\closeout\reffile

\centerline{\bf SPECTRAL SYNTHESIS AND TOPOLOGIES ON IDEAL SPACES} 
\centerline{\bf FOR BANACH $^*$-ALGEBRAS} 
\bigskip 
\centerline{\bf J. F. Feinstein, E. Kaniuth, and D. W. B. Somerset} 
\bigskip 
\bigskip 
\noindent{\bf Abstract} This paper continues the study of spectral synthesis 
and the topologies $\tau_{\infty}$ and $\tau_r$ on 
the ideal space of a Banach algebra, concentrating on the class of Banach 
$^*$-algebras, and in particular on $L^1$-group 
algebras. It is shown that if a group $G$ is a finite 
extension of an abelian group then $\tau_r$ is Hausdorff on the ideal space of 
$L^1(G)$ if and only if $L^1(G)$ has spectral 
synthesis, which in turn is equivalent to $G$ being compact. The result is 
applied to nilpotent groups, [FD]$^-$-groups, and 
Moore groups. An example is given of a non-compact, 
non-abelian group $G$ for which 
$L^1(G)$ has spectral synthesis. It is also shown that 
if $G$ is a non-discrete group then $\tau_r$ is not Hausdorff on the ideal 
lattice of the Fourier algebra $A(G)$. 
\bigskip 
\noindent{{\bf 1991 Maths Subject Classification} 46H10, 46K05, 43A45, 22D15} 
\bigskip 
\bigskip 
\noindent {\bf 1. Introduction} 
\bigskip 
\noindent The goal of ideal theory for Banach algebras is, on the one hand to 
describe the set $Id(A)$ of closed two-sided 
ideals of a Banach algebra $A$, and on the other hand to use knowledge of the 
ideal structure to obtain information about the 
algebra itself. This usually involves representing the algebra as a bundle of 
quotient algebras over a topological base-space. 
Standard examples include the Gelfand theory for commutative Banach algebras 
and the theory of continuous bundles of 
C$^*$-algebras. Another familiar example is the representation of a 
C$^*$-algebra as a bundle of C$^*$-algebras over its 
primitive ideal space with the hull-kernel topology. 

In trying to describe the closed ideals of a Banach algebra $A$, the first 
question is whether it is sufficient to know about 
the primitive ideals of $A$. This is the question of `spectral synthesis', 
originally studied for commutative Banach algebras, 
then for $L^1$-group algebras and other Banach $^*$-algebras, and more 
recently also for the Haagerup tensor product of 
C$^*$-algebras [\the\ASS], [\the\fes]. On the other hand, the representing of 
a Banach algebra as a bundle of quotient algebras 
involves the study of topologies on the space $Id(A)$, and its subsets. 
This is an area which has received considerable attention in recent years, 
with the work of Archbold [\the\RJA], Beckhoff [\the\Be], [\the\Bec], 
[\the\Beck], and others [\the\AKSS], [\the\HKMS], 
[\the\Id]. 

The two main topologies that have been introduced are $\tau_{\infty}$ and 
$\tau_r$. The first, $\tau_{\infty}$, was defined by 
Beckhoff [\the\Be], using the various continuous norms and seminorms which the 
algebra can carry. It is always compact on 
$Id(A)$, but seldom Hausdorff. In fact, for a commutative Banach algebra $A$, 
$\tau_{\infty}$ is Hausdorff on $Id(A)$ if and 
only if $A$ has spectral synthesis [\the\fes; 1.8]. Since the topology 
$\tau_{\infty}$ is also Hausdorff for C$^*$-algebras, 
and since C$^*$-algebras also have a form of spectral synthesis, the third 
author was led to introduce a definition of spectral 
synthesis for a general Banach algebra $A$, and to investigate the relation 
between this and $\tau_{\infty}$ [\the\Syn], 
[\the\fes]. It is known that if $A$ is a separable, unital PI-Banach algebra 
then $A$ has spectral synthesis if and only if 
$\tau_{\infty}$ is Hausdorff on $Id(A)$, and that for a general Banach algebra 
$A$ the Hausdorffness of $\tau_{\infty}$ on 
$Id(A)$ implies a weak form of spectral synthesis. Conversely, if $A$ is 
separable and has a strong form of spectral synthesis 
then $\tau_{\infty}$ is a $T_1$-topology. Thus spectral synthesis and the 
Hausdorffness of $\tau_{\infty}$ seem to be closely 
related, and possibly identical. 

The second topology that has been introduced is $\tau_r$ [\the\Id]. This is 
also compact on $Id(A)$, and is Hausdorff whenever 
$\tau_{\infty}$ is [\the\Id; 3.1.1]. It was shown in [\the\Id; 2.11] that if 
there is a compact Hausdorff topology on a 
subspace of $Id(A)$, which is related to the quotient norms in a useful way, 
then that topology necessarily coincides with the 
restriction of $\tau_r$. Thus, for instance, $\tau_r$ coincides with the 
Gelfand topology on the maximal ideal space of a 
unital commutative Banach algebra. There are a number of cases when $\tau_r$ 
is Hausdorff but $\tau_{\infty}$ is not, e.g. for 
TAF-algebras [\the\Id] and for the Banach algebra $C^1[0,1]$ [\the\FS]. For 
uniform algebras, however, it turned out that 
$\tau_r$ is Hausdorff if and only if $\tau_{\infty}$ is Hausdorff [\the\FS]. 
Thus for uniform algebras without spectral 
synthesis, such as the disc algebra, there is no useful compact Hausdorff 
topology on the space of closed ideals. 

The purpose of this paper is to continue the study of spectral synthesis and 
the topologies $\tau_{\infty}$ and $\tau_r$, 
concentrating on the class of Banach $^*$-algebras. We are interested in the 
questions of whether spectral synthesis is 
equivalent to the Hausdorffness of $\tau_{\infty}$ for these algebras, and 
whether $\tau_r$ can be Hausdorff when spectral 
synthesis fails. In Section 2 we employ general techniques of Banach 
$^*$-algebra theory. The main result of the section is 
that spectral synthesis is equivalent to the Hausdorffness of $\tau_{\infty}$ 
for a class of Banach $^*$-algebras which 
includes the $L^1$-algebras of [FC]$^-$-groups. In Section 3 we employ 
group-theoretic techniques. We give an example of a 
non-compact, non-abelian 
group $G$ for which $L^1(G)$ has spectral synthesis, but we show 
that for nilpotent groups, [FD]$^-$-groups, and 
Moore groups, spectral synthesis for $L^1(G)$ is equivalent to the compactness 
of $G$, and furthermore that if $G$ fails to be 
compact then $\tau_r$ fails to be Hausdorff on $Id(L^1(G))$. We also show that 
if $G$ is a non-discrete group then $\tau_r$ is 
not Hausdorff on the ideal lattice of the Fourier algebra $A(G)$. 

We now give the definitions of the various topologies on $Id(A)$, starting 
with the lower topology $\tau_w$. Let $A$ be a 
Banach algebra. A subbase for $\tau_w$ on $Id(A)$ is given by the sets $\{ 
I\in Id(A): I \not\supseteq J\}$ as $J$ varies 
through the elements of $Id(A)$. Thus the restriction of $\tau_w$ to the set 
of closed prime ideals is simply the hull-kernel 
topology. Next we define $\tau_{\infty}$. For each $k\in {\bf N}$, let 
$S_k=S_k(A)$ denote the set of seminorms (`seminorm' 
means `algebra seminorm' in this paper) $\rho$ on $A$ satisfying $\rho(a)\le 
k\Vert a\Vert$ for all $a\in A$. Then $S_k$ is a 
compact, Hausdorff space [\the\Be]. We say that $\rho\ge\sigma$, for 
$\rho,\sigma\in S_k$, if $\rho (a)\le\sigma (a)$ for all 
$a\in A$. The point of this upside-down definition is that if $\rho\ge\sigma$ 
then $\ker\rho\supseteq\ker\sigma$. Clearly if 
$\rho,\sigma\in S_k$ the seminorm $\rho\wedge\sigma$ defined by 
$(\rho\wedge\sigma)(a)=\max\, \{\rho(a),\sigma(a)\}$ is the 
greatest seminorm less than both $\rho$ and $\sigma$ in the order structure. 
Thus $S_k$ is a lattice. The topology 
$\tau_{\infty}$ is defined on $Id(A)$ as follows [\the\Be]: for each $k$ let 
$\kappa_k:S_k\to Id(A)$ be the map 
$\kappa_k(\rho)=\ker\rho$, and let $\tau_k$ be the quotient topology of 
$\kappa_k$ on $Id(A)$. Then 
$\tau_{\infty}=\bigcap_k\tau_k$. Clearly each $\tau_k$ is compact, so 
$\tau_{\infty}$ is compact. 

Next we define the topology $\tau_r$, which is the join of two weaker 
topologies. The first is easily defined: $\tau_u$ is the weakest topology on 
$Id(A)$ for which all the norm functions $I\mapsto\Vert a+I\Vert$ $(a\in A,\ 
I\in Id(A))$ are upper semi-continuous. The other topology $\tau_n$ can be 
described in various different ways, but none is particularly easy to work 
with. A net $( I_{\alpha})$ in $Id(A)$ is said to have the {\sl normality 
property} with respect to $I\in Id(A)$ if $a\notin I$ implies that 
$\lim\inf\Vert a+I_{\alpha}\Vert >0$. Let $\tau_n$ be the topology whose 
closed sets $N$ have the property that if $(I_{\alpha})$ is a net in $N$ with 
the normality property relative to $I\in Id(A)$ 
then $I\in N$. It follows that if $(I_{\alpha})$ is a net in $Id(A)$ having 
the normality property relative to 
$I\in Id(A)$ then $I_{\alpha}\to I$ $(\tau_n)$. Any topology for which 
convergent nets have the normality property with respect to each of their 
limits (such a topology is said to have the {\sl normality property}) is 
necessarily stronger than $\tau_n$, but $\tau_n$ itself need not have the 
normality property. Indeed the following is true. Let 
$\tau_r$ be the topology on $Id(A)$ generated by $\tau_u$ and $\tau_n$. Then 
$\tau_r$ is always compact [\the\Id; 2.3], and 
$\tau_r$ is Hausdorff if and only if $\tau_n$ has the normality property 
[\the\Id; 2.12]. It is a useful fact that for $I\in 
Id(A)$, $Id(A/I)$ is $\tau_{\infty}$-- $\tau_{\infty}$ and $\tau_r$-- $\tau_r$ 
homeomorphic to the subset $\{ J\in 
Id(A):J\supseteq I\}$ of $Id(A)$ [\the\Be; Prop. 5], [\the\Id; 2.9].

The following simple lemma is taken from [\the\FS; 0.1]. 
\bigskip 
\noindent {\bf Lemma 1.1} {\sl Let $A$ be a Banach algebra. Let $(I_{\alpha})$ 
be a net in $Id(A)$, either 
decreasing or increasing, and correspondingly either set $I=\bigcap 
I_{\alpha}$ 
or $I=\overline{\bigcup I_{\alpha}}$. Then $I_{\alpha}\to I$ $(\tau_r)$.} 
\bigskip 
\noindent Now let $A$ be a Banach algebra and let $Prim(A)$ be the space of 
primitive ideals of $A$ (i.e. the kernels of 
algebraicially irreducible representations of $A$) equipped with the 
hull-kernel topology. Let $Prime(A)$ be the space of 
proper closed prime ideals of $A$, and let $Prim^s(A)$ be the space of 
semisimple prime ideals of $A$ (such ideals are 
automatically closed), both spaces also being equipped with the hull-kernel 
topology. The notation, and the importance of 
$Prim^s(A)$, is explained in [\the\fes]. 

The paper [\the\Syn] contained a definition of `spectral synthesis', but 
unfortunately that definition was slightly too 
restrictive, and was replaced in [\the\fes] by the following definition. 
\bigskip 
\noindent {\bf Definition of spectral synthesis} A Banach algebra $A$ has {\sl 
spectral synthesis} if it has the following 
properties: 

(i) $Prim^s(A)$ is locally compact. 

(ii) $\tau_w$ has the normality property on $Prim(A)$ (or equivalently on 
$Prim^s(A)$). 

(iii) $Id(A)$ is isomorphic to the lattice of open subsets of $Prim(A)$, under 
the correspondence $I\leftrightarrow \{ P\in 
Prim(A):P\not\supseteq I\}$. 
Equivalently, every proper, closed ideal of $A$ is semisimple. 

\bigskip 
\noindent {\bf Remarks} (a) It was noted in [\the\fes], see [\the\Syn; 1.1], 
that if $A$ has spectral synthesis then $\tau_w$ 
has the normality property on the whole of $Id(A)$. 

(b) The definition just given coincides with the standard one for the class of 
commutative Banach algebras [\the\fes; 1.7]. 
Recall that a (possibly non-unital) commutative Banach algebra $A$ has 
spectral synthesis (usual definition) if the map 
$I\mapsto \{P\in Prim(A):P\supseteq I\}$ sets up a 1--1 correspondence between 
closed ideals of $A$ and Gelfand closed subsets 
of $Prim(A)$. This is equivalent to requiring that the hull-kernel and Gelfand 
topologies coincide on $Prim(A)$, and that every 
closed ideal of $A$ is semisimple. 
\bigskip 
\bigskip 
\noindent{\bf 2. Banach $^*$-algebras} 
\bigskip 
\noindent In this section we consider spectral synthesis and the topologies 
$\tau_{\infty}$ and $\tau_r$ within the class of 
Banach $^*$-algebras. We show that our definition of spectral synthesis 
coincides with the usual definition for a large 
subclass which probably contains all the cases of interest. For a smaller 
class, which contains the $L^1$-algebras of 
[FC]$^-$-groups, we are able to show that spectral synthesis is equivalent to 
the Hausdorffness of $\tau_{\infty}$. 

Let $A$ be a $^*$-semisimple Banach $^*$-algebra with C$^*$-envelope $C^*(A)$ 
and self-adjoint part $A_{sa}$. Let $Prim_*(A)$ 
be the set of kernels of topologically irreducible $*$-representations of $A$ 
on Hilbert space, with the hull-kernel topology 
(such representations are automatically continuous [\the\BD; p.196]). The map 
$P\mapsto P\cap A$ $(P\in Prim(C^*(A)))$ maps 
$Prim(C^*(A))$ continuously onto $Prim_*(A)$, but is not a homeomorphism in 
general. If it is a homeomorphism then $A$ is said 
to be {\sl $^*$-regular}. A Banach $^*$-algebra $A$ is said to be {\sl 
hermitian} (or {\sl symmetric}) if every self-adjoint 
element of $A$ has real spectrum. This has the implication that every 
primitive ideal of $A$ is the kernel of a topologically 
irreducible $*$-representation on a Hilbert space, or in other words, that 
$Prim(A)\subseteq Prim_*(A)$, see [\the\Lept; 
pp.50-51]. 

Let $A$ be a $^*$-semisimple Banach $^*$-algebra. The `usual definition of 
spectral synthesis' is that $A$ has spectral 
synthesis if each closed subset of $Prim_*(A)$ is the hull of a unique closed 
ideal of $A$. We begin by showing that for a 
large class of Banach $^*$-algebras, which probably contains all the relevant 
examples, this `usual definition' is equivalent 
to our definition. Our method requires $A$ to be hermitian and $^*$-regular, 
but spectral synthesis (in either sense) is a very 
strong property, and it seems unlikely that a Banach $^*$-algebra could have 
spectral synthesis but fail to be hermitian and 
$^*$-regular. 

\bigskip 
\noindent {\bf Theorem 2.1} {\sl Let $A$ be a hermitian, $^*$-regular, 
$^*$-semisimple Banach $^*$-algebra. Then $A$ has 
spectral synthesis in the sense of this paper if and only if $A$ has spectral 
synthesis in the usual sense for $^*$-semisimple 
Banach $^*$-algebras.} 
\bigskip 
\noindent {\bf Proof.} As we noted above, we have $Prim(A)\subseteq 
Prim_*(A)$, since $A$ is hermitian. We also have that 
$Prim_*(A)\subseteq Prim^s(A)$ since every $P\in Prim_*(A)$ is prime and 
semisimple. 

Now if $A$ has spectral synthesis in our sense then, by condition (iii), the 
elements of $Prim(A)$ separate the closed ideals 
of $A$. Since $Prim(A)\subseteq Prim_*(A)$, it follows that the elements of 
$Prim_*(A)$ separate the closed ideals of $A$, and 
thus that $A$ has spectral synthesis in the usual sense. Conversely, if $A$ 
has spectral synthesis in the usual sense then the 
elements of $Prim_*(A)$ separate the closed ideals of $A$. Since each element 
of $Prim_*(A)$ is semisimple, it follows that the 
primitive ideals of $A$ separate the closed ideals of $A$, and thus condition 
(iii) of our version of spectral synthesis holds. 
It remains to show that conditions (i) and (ii) of our definition of spectral 
synthesis are automatically satisfied. 

Since $A$ is $^*$-regular, $Prim_*(A)$ is locally compact, which implies that 
$Prim^s(A)$ is locally compact, by Remark (a) 
after [\the\fes; Definition 1.2]. Thus property (i) of spectral synthesis 
holds. 

For property (ii), first observe that if the normality property fails on 
$Prim(A)$ for an element $a\in A$, then it also fails 
for $a^*a$. For suppose that $(P_{\alpha})$ is a net in $Prim(A)$ converging 
to $P\in Prim(A)$, and that $a\notin P$, but that 
$\lim_{\alpha}\Vert a+P_{\alpha}\Vert=0$. Then $a^*a\notin P$, because $P$ is 
the kernel of a topologically irreducible 
$*$-representation of $A$ on Hilbert space, but $\lim_{\alpha}\Vert 
a^*a+P_{\alpha}\Vert=0$. Thus it is enough to establish 
that the normality property holds for self-adjoint elements. Let $a\in A_{sa}$ 
and let $P\in Prim_*(A)$. Then the quotient norm 
of $a+P$ in $A/P$ dominates the spectral radius of $a+P$ in $A/P$. Let $C$ be 
the completion of $A/P$ in the C$^*$-norm from 
the corresponding topologically irreducible $^*$-representation of $A/P$. Then 
the spectral radius of the canonical image of 
$a+P$ in $C$ is less than or equal to the spectral radius of $a+P$ in $A/P$. 
But the C$^*$-norm of $a+P$ is equal to the 
spectral radius in $C$. Thus we have shown that, for self-adjoint elements, 
the quotient norm dominates the corresponding 
C$^*$-norm, for each $P\in Prim_*(A)$. On the other hand, the norm functions 
are lower semicontinuous on $Prim_*(A)$ for the 
C$^*$-norms, since $Prim_*(A)\cong Prim(A)$, and it is straightforward to show 
from this that the normality property holds for 
the quotient norms on $Prim_*(A)$. This establishes (ii). Q.E.D. 
\bigskip 
\noindent If $G$ is a locally compact group $G$, then $A=L^1(G)$ is a 
$^*$-semisimple Banach $^*$-algebra, and $C^*(A)$ is the 
(full) group C$^*$-algebra $C^*(G)$ of $G$. The class of locally compact 
groups for which $L^1(G)$ is hermitian and 
$^*$-regular includes all nilpotent groups and all connected groups of 
polynomial growth [\the\Lu], [\the\BLSV]. It also 
includes all groups in [FC]$^-$ [\the\HKK], where a group belongs to [FC]$^-$ 
provided that each conjugacy class has compact 
closure. 

For a locally compact abelian group $G$, $L^1(G)$ has spectral synthesis if 
and only if $G$ is compact. This classical theorem 
is chiefly due to Malliavin. 
In the next section we exhibit a non-compact, non-abelian
group $G$ for which $L^1(G)$ has 
spectral synthesis. 
\bigskip 
\noindent It is not difficult to see that an algebra $A$ in the class of 
hermitian, $^*$-regular, $^*$-semisimple Banach 
$^*$-algebras has spectral synthesis if and only if every closed ideal of $A$ 
is semisimple. Thus from [\the\fes; 1.15] we have 
the following useful extension result. 
\bigskip 
\noindent {\bf Proposition 2.2} {\sl Let $A$ be a hermitian, $^*$-regular, 
$^*$-semisimple Banach $^*$-algebra, and let $J\in 
Id(A)$. If both $J$ and $A/J$ have spectral synthesis then $A$ has spectral 
synthesis.} 
\bigskip 
\noindent {\bf Proof.} It follows from the remarks after [\the\fes; 1.15] that 
every closed ideal of $A$ is semisimple. Hence $A$ has spectral synthesis, as 
we have just observed. Q.E.D. 
\bigskip 
\noindent We are interested in showing that spectral synthesis is equivalent 
to the Hausdorffness of $\tau_{\infty}$ on 
$Id(A)$. The next result establishes one direction of this, for a particular 
class of Banach $^*$-algebras. 
\bigskip 
\noindent {\bf Theorem 2.3} {\sl Let $A$ be a $^*$-semisimple Banach 
$^*$-algebra, and suppose that $Prim(C^*(A))$ is 
Hausdorff. If $\tau_{\infty}$ is Hausdorff on $Id(A)$ then $A$ has spectral 
synthesis (in the sense of this paper).} 
\bigskip 
\noindent {\bf Proof.} By [\the\fes; 1.11] it is enough to show that every 
proper, closed, prime ideal of $A$ belongs to 
$Prim_*(A)\subseteq Prim^s(A)$. Thus let $P$ be a proper, closed, prime ideal 
of $A$. Since $\bigcap \{ Q:Q\in Prim_*(A)\}=\{ 
0\}$, $Prim_*(A)$ is $\tau_w$-dense in $Prime(A)$. Thus there is a net 
$(P_{\alpha})$ in $Prim_*(A)$ converging to $P$ 
$(\tau_w)$. Let $(q_{\alpha})$ be the corresponding net of quotient norms in 
the compact space $S_1(A)$. By passing to a 
subnet, if necessary, we may assume that $(q_{\alpha})$ converges to some 
seminorm $q$, say, in $S_1(A)$. Set $Q=\ker q$. Then 
$P_{\alpha}\to Q$ $(\tau_{\infty})$, and the normality property for $\tau_w$ 
[\the\Syn; 2.5] implies that $Q\subseteq P$. Set 
$B=C^*(A)$. For each $\alpha$, let $\widetilde P_{\alpha} \in Prim(B)$ such 
that $\widetilde P_{\alpha}\cap B=P_{\alpha}$. By 
the $\tau_{\infty}$-compactness of $Id(B)$ we may assume that the net 
$(\widetilde P_{\alpha})$ is $\tau_{\infty}$-convergent 
in $Id(B)$, with limit $\widetilde R$ say. Since $Prim(B)$ is Hausdorff, the 
set $Prim(B)\cup \{B\}$ is $\tau_{\infty}$-closed 
in $Id(B)$, by [\the\Be; Proposition 8] and [\the\RJA; 3.3(b), 4.3(b)]. Hence 
either $\widetilde R=B$, or $\widetilde R\in 
Prim(B)$. But since $\widetilde P_{\alpha} \cap A=P_{\alpha}\to \widetilde R 
\cap A$ $(\tau_{\infty})$ (the restrictions to A 
of the quotient C$^*$-seminorms converge), we have that $\widetilde R \cap 
A=Q$, by the Hausdorffness of $\tau_{\infty}$ on 
$Id(A)$. Hence $\widetilde R\ne B$ so $\widetilde R\in Prim(B)$. Thus $P\in 
Prim_*(A)$, as required. Q.E.D. 
\bigskip 
\noindent It was shown in [\the\Kant] that if $G$ is an [FC]$^-$-group then 
$Prim(C^*(G))$ is Hausdorff. 
\bigskip 

Now we work in the other direction, trying to show that spectral synthesis 
implies that $\tau_{\infty}$ is Hausdorff. The 
following definition is taken from [\the\fes]. 

Let $A$ be a Banach algebra, and suppose that there is a continuous norm 
$\gamma$ on $A$ such that $B$, the $\gamma$-completion 
of $A$, is a C$^*$-algebra (in this paper $A$ will be a Banach $^*$-algebra 
with $B=C^*(A)$, so $\gamma$ will be the maximal 
C$^*$-seminorm on $A$). Extending the definition from the Haagerup tensor 
product of C$^*$-algebras [\the\ASS; \S6], we refer 
to those closed ideals in $A$ of the form $J\cap A$ $(J\in Id(B))$ as {\sl 
upper ideals}. The set of upper ideals is denoted 
$Id^u(A)$. Note that if $I$ is an upper ideal of $A$ then in fact $I=J\cap A$ 
where $J$ is the closure of $I$ in $B$. 

\bigskip 
\noindent {\bf Definition} [\the\fes; 2.4] Let $A$ be a Banach algebra. We 
shall say that $A$ has {\sl property (P)} if $A$ 
satisfies the following conditions: 

(a) there is a continuous norm $\gamma$ on $A$ such that $B$, the 
$\gamma$-completion of $A$, is a C$^*$-algebra; 

(b) every primitive ideal of $A$ is an upper ideal, i.e. $Prim(A)\subseteq 
Id^u(A)$; 

(c) there is a subset $R$ of $A\cap B_{sa}$ such that each $a\in R$ is 
contained in a completely regular, commutative Banach 
$^*$-subalgebra $A_a$ of $A$ (where the norm and the involution on $A_a$ are 
those induced by $B$), and such that if $I\in 
Id^u(A)$ and $J\in Id(A)$ with $J\not\subseteq I$ then there exists $a\in R$ 
such that $A_a\cap J\not\subseteq I$. 
\bigskip 
\noindent Recall that a $^*$-semisimple Banach $^*$-algebra $A$ is said to be 
{\sl locally regular} if there is a dense subset $R$ of $A_{sa}$ such that 
each $a\in R$ generates a completely regular 
commutative Banach $^*$-subalgebra $A_a$ [\the\Bar; \S 4]. A locally regular, 
$^*$-semisimple Banach $^*$-algebra is 
automatically $^*$-regular [\the\Bar; Theorem 4.3]. It is shown in [\the\Bar; 
Theorem 4.1] that if a locally compact group $G$ 
has polynomial growth then $L^1(G)$ is locally regular. Furthermore, if 
$\omega$ is a polynomial weight and either $G$ is 
compactly generated and of polynomial growth [\the\Bar], or $G\in$[SIN] 
[\the\Run], then the Beurling algebra $L^1(G,\omega)$ 
is also locally regular. 
\bigskip 
\noindent {\bf Theorem 2.4} {\sl Let $A$ be a hermitian, locally regular, 
$^*$-semisimple Banach $^*$-algebra. Then $A$ has 
property (P). Hence every upper ideal of $Id(A)$ is $\tau_{\infty}$-closed in 
$Id(A)$, and if $A$ has spectral synthesis then 
$\tau_{\infty}$ is Hausdorff on $Id(A)$. 
If $Prim(C^*(A))$ is Hausdorff then $Id(A)$ is $\tau_{\infty}$-Hausdorff if 
and only if $A$ has spectral synthesis.} 
\bigskip 
\noindent {\bf Proof.} First we show that $A$ has property (P). Condition (a) 
is immediate (see [\the\BD; p.223] if need be). 
Condition (b) follows from the fact that $A$ is hermitian, see [\the\Lept; 
pp.50-51]. 
It is shown in [\the\HKK ; Lemma 1.2] that if $A$ is a hermitian, locally 
regular, $^*$-semisimple Banach $^*$-algebra then for 
each closed subset $E$ of $Prim(A)$ there is a smallest ideal $J(E)$ of $A$ 
with hull equal to $E$. Furthermore, $J(E)$ is 
generated by elements from the set $\bigcup\{ A_a:a\in R\}$. Let $I$ be an 
upper ideal of $A$ and let $J$ be a closed ideal of 
$A$ not contained in $I$. Set $E=\{ P\in Prim(A): P\supseteq I\}$ and $F=\{ 
P\in Prim(A): P\supseteq J\}$. Then $E\not\subseteq 
F$, so [\the\HKK; Lemma 1.2] shows that $J(F)\not\subseteq I$. Hence there 
exists $a\in R$ such that $A_a\cap J(F)\not\subseteq 
I$. Since $J(F)\subseteq J$, we have established that condition (c) holds. 

It now follows from [\the\fes; Proposition 2.6] that every upper ideal of 
$Id(A)$ is $\tau_{\infty}$-closed in $Id(A)$, and if 
$A$ has spectral synthesis then $\tau_{\infty}$ is Hausdorff on $Id(A)$. If, 
furthermore, $Prim(C^*(A))$ is Hausdorff then it 
follows from [\the\fes; Proposition 2.6] and Theorem 2.3. 
that $Id(A)$ is $\tau_{\infty}$-Hausdorff if and only if $A$ has spectral 
synthesis. Q.E.D. 
\bigskip 
\noindent {\bf Corollary 2.5} {\sl Let $G$ be an [FC]$^-$-group. Then 
$\tau_{\infty}$ is Hausdorff on $Id(L^1(G))$ if and only 
if $L^1(G)$ has spectral synthesis.} 
\bigskip 
\noindent {\bf Proof.} $L^1(G)$ is hermitian and locally regular, and 
$Prim(C^*(G))$ is Hausdorff [\the\Kant], so the result 
follows from Theorem 2.4. Q.E.D. 
\bigskip 
\noindent In the same way it also follows from Theorem 2.4 that if $G$ is a 
connected group of polynomial growth then every 
upper ideal of $L^1(G)$ is $\tau_{\infty}$-closed in $Id(L^1(G))$. 

It is unknown, for the class of [FC]$^-$-groups, whether spectral synthesis 
for $L^1(G)$ is equivalent to the compactness of 
$G$. Indeed very little is known about spectral synthesis for [FC]$^-$-groups. 
It is not even known whether singletons in the 
primitive ideal space are sets of synthesis. 
\bigskip 
\bigskip 
\noindent {\bf 3. L$^1$-group algebras and Fourier algebras} 
\bigskip 
\noindent In this section we employ group-theoretic techniques, concentrating 
mainly on $L^1$-group algebras. We begin with an 
example of a non-compact, non-abelian 
group $G$ for which $L^1(G)$ has spectral synthesis. 
Next we show that if $G$ is a finite extension of 
an abelian group then the topology $\tau_r$ is Hausdorff on $Id(L^1(G))$ if 
and only if $L^1(G)$ has spectral synthesis, which 
occurs if and only if $G$ is compact. General results allow one to apply this 
to the classes of nilpotent groups, 
[FD]$^-$-groups, and Moore groups (where [FD]$^-$ denotes the class of locally 
compact groups for which the commutator subgroup 
has compact closure, and a locally compact group $G$ is a Moore group if every 
irreducible unitary representation of $G$ is 
finite dimensional). Finally we show that if $G$ is a non-discrete group then 
$\tau_r$ fails to be Hausdorff on the ideal space 
of the Fourier algebra $A(G)$. 
\bigskip 

\noindent {\bf Example 3.1} {\sl A non-compact group with spectral synthesis}. 
Let $p$ be a prime and let $N$ be the field of $p$-adic numbers. Let $K$ be 
the subset of elements of $N$ of valuation $1$. 
Then $K$ is a compact subgroup under multiplication. Let $G=K\propto N$, where 
$K$ acts on the additive group $N$ by 
multiplication. The group $G$ is often referred to as Fell's non-compact group 
with countable dual. 

The dual group $\widehat N$ is canonically isomorphic to $N$. More precisely, 
there is a character $\chi$ of $N$ such that 
$\chi (x)=1$ if and only if $x$ is a $p$-adic integer, and then the mapping 
$y\mapsto \chi_y$, where $\chi_y(x)=\chi(xy)$, is a 
topological isomorphism between $N$ and $\widehat N$. The irreducible 
representations of $G$ are easy to determine using 
Mackey's theory. For $y\in N_j$, the $K$-orbit of $\chi_y$ is equal to $\{ 
\chi_t:t\in N_j\}$. Let $\pi_j={\rm ind}^G_N 
\chi_y$, for $y\in N_j$. Then $\widehat G=\widehat K \cup\{\pi_j:j\in {\bf 
Z}\}$. The topology of $\hat G$ has been described 
in [\the\Ba; 4.6]. Both $\widehat K$ and $\widehat G\backslash \widehat K$ are 
discrete, and a sequence $(\pi_{j_k})_k$ 
converges to some (and hence all) $\sigma\in\widehat K$ if and only if $j_k\to 
-\infty$. 

Now let $E$ be a closed subset of $\widehat G$, and let $J_E=\{ j\in {\bf 
Z}:\pi_j\in E\}$ and $$F=\{ 0\}\cup\bigcup_{j\in 
J_E}\{\chi_y:y\in N_j\}\subseteq \widehat N.$$ 
Then $F$ is closed and $G$-invariant, and since the $N_j$ are open and closed 
in $\widehat N$, the boundary of $F$ is contained 
in the singleton $\{ 0\}$. Thus $F$ is a spectral set for $L^1(N)$. The 
projection theorem for spectral sets [\the\HaLu] shows 
(in the notation of [\the\HaLu]) that $$h(e_N(k(F)))=\widehat K\cup E$$ is a 
spectral set for $L^1(G)$. Since $\widehat K$ is 
discrete, $E$ is open (and closed) in $\widehat K \cup E$, so since $L^1(G)$ 
has the Wiener property [\the\Lu], it follows that 
$E$ is a spectral set for $L^1(G)$, see [\the\HKK; Remark 1.3]. 

\bigskip 
\noindent This example is hermitian by [\the\LeP] and $^*$-regular, so the two 
possible versions of spectral synthesis 
coincide, by Theorem 2.1. 
The group $G$ also has polynomial growth, so $L^1(G)$ is locally regular. 
Hence the topologies $\tau_{\infty}$ and $\tau_r$ are 
Hausdorff on $Id(L^1(G))$, by Theorem 2.4 and [\the\Id; 3.1.1]. 

Another strategy for showing that $L^1(G)$ above has spectral synthesis might 
be to proceed as follows. Let $J$ be the ideal of 
$A=L^1(G)$ given by $J=\ker \widehat K$. Then $J$ is a semisimple Banach 
algebra with discrete primitive ideal space. If it 
could be shown directly that $J$ has spectral synthesis then since $A/J\cong 
L^1(K)$ also has spectral synthesis (because $K$ 
is compact), it would follow from Proposition 2.2 that $A$ has spectral 
synthesis. 

\bigskip 
\noindent We now go on to show that the failure of spectral synthesis implies 
that the topology $\tau_r$ is non-Hausdorff, at 
least for a large class of groups. We are grateful to Colin Graham for 
informing us about the following proposition, which he 
can prove by tensor methods. Here we provide an alternative proof, in keeping 
with the methods of this paper. 
\bigskip 
\noindent{\bf Proposition 3.2} {\sl Let $G$ be a non-compact locally compact 
abelian group. Then every non-empty open subset of 
$\widehat G$ contains a closed subset which is non-spectral.} 
\bigskip 
\noindent {\bf Proof.} Notice first that for $\alpha\in\widehat G$ and $I\in 
Id(L^1(G))$, the mapping $f\mapsto \alpha f$ is an 
automorphism of $L^1(G)$ and $h(\alpha I)=\alpha^{-1}h(I)$. It is enough 
therefore to show that every open neighbourhood $V$ of 
$1$ contains a closed non-spectral set. 

Since the topology on $\widehat G$ is the topology of uniform convergence on 
compact subsets of $G$, we can assume that $V$ is 
of the form 
$$V=\{ \alpha\in\widehat G: |\alpha(x)-1|<\epsilon\hbox{ for all }x\in C\},$$ 
where 
$\epsilon>0$ and $C$ is a compact subset of $G$. Let $H$ be any 
compactly-generated (open) subgroup of $G$ containing $C$, and 
let $\phi:\widehat G\to \widehat H$ denote the restriction map 
$\alpha\mapsto\alpha|_H$. Then $$ 
V=\phi^{-1}\left(\{\gamma\in \widehat H:|\gamma(x)-1|<1\hbox{ for all }x\in 
C\}\right).$$ 
In particular, when $H$ is compact, $\widehat{G/H}$ is open in $\widehat G$ 
and $\widehat{G/H}\subseteq V$. Since, in this 
case, $G/H$ is non-compact, $\widehat{G/H}$ has a closed subset which is 
non-spectral for $L^1(G/H)$, and hence is non-spectral 
for $L^1(G)$ by the injection theorem for spectral sets. Thus we can 
henceforth assume that $H$ is non-compact. 

By the projection theorem for spectral sets, for any closed subset $F$ of 
$\widehat H$, $F$ is a spectral set for $L^1(H)$ if 
and only if $\phi^{-1}(F)$ is a spectral set for $L^1(G)$. Therefore it 
suffices to treat the case when $G$ is non-compact and 
compactly generated. By the structure theorem, $$G={\bf R}^m\times {\bf 
Z}^n\times K,$$ where $K$ is a compact group and 
$m+n\ge 1$. When $m\ge 1$, pass to $H=G/{\bf R}^{m-1}\times {\bf Z}^n\times 
K={\bf R} $, and when $m=0$ pass to $H=G/{\bf 
Z}^{n-1}\times K={\bf Z}$. Suppose that we know that every neighbourhood of 
$1$ in $\widehat H$ contains a closed non-spectral 
set. Then applying the injection theorem for spectral sets once more, the same 
follows for $\widehat G$. Thus we are reduced to 
the two cases $G={\bf R}$ and $G={\bf Z}$. 

Consider $G={\bf R}$ first. It is well-known that $\widehat {\bf R} ={\bf R}$ 
contains a compact non-spectral set, although we have not been able to find
a specific reference for this fact (in fact such a set can be formed as a 
finite union of translates by integers of a suitable compact subset $E$ of 
$[0,1]$, where $E$ is such that the set $\{\exp(2\pi i s): s\in E\}$ is
a non-spectral set for $L^1(\bf Z)$). Thus it is enough to show 
that if $E\subseteq {\bf R}$ is non-spectral then so is $sE$ for every $s>0$. 
Now it is easily verified that the mapping 
$\theta:L^1({\bf R})\to L^1({\bf R})$ given by 
$$(\theta f)(x)=\frac1{s} f\left(\frac{x}{s}\right),\ \ \ x\in{\bf R},$$ 
is an (isometric) isomorphism of $L^1({\bf R})$ and satisfies 
$$\widehat{\theta f}(y)=\widehat f (sy),\ \ \ y\in{\bf R}.$$ 
It follows that $sE$ is non-spectral whenever $E$ is. 

Finally, let $G={\bf Z}$. Identifying $\widehat{\bf R}$ with ${\bf R}$ and 
$\widehat{\bf Z}$ with ${\bf T}$, the restriction 
map $p:\widehat{\bf R}\to \widehat{\bf Z}$ is given by $p(y)=e^{2\pi iy}$. Now 
if $V$ is any neighbourhood of $1$ in ${\bf T}$, 
choose $0<\delta<\frac12$ such that $p( (-\delta, \delta))\subseteq V$. Let 
$E$ be a closed subset of ${\bf R}$ such that $E$ 
is non-spectral for $L^1({\bf R})$ and $E\subseteq (-\delta, \delta)$, and let 
$F=p(E)\subseteq V$. Then $E$ is open (and 
closed) in $p^{-1}(F)=\bigcup_{m\in{\bf Z}} (m+E)$. Since a clopen subset of a 
spectral set is itself spectral, it follows that 
$p^{-1}(F)$ is not spectral. Hence the projection theorem implies that $F$ is 
non-spectral. Q.E.D. 
\bigskip 
\noindent In order to exploit the existence of non-spectral sets, we need some 
information about quotient norms. The following 
definition is useful. 
  
Let $A$ be a completely regular, natural Banach function algebra on its 
maximal ideal space $Max(A)$. Recall that a Gelfand 
compact subset $X$ of $Max(A)$ is a {\sl Helson set} if $A|_X=C(X)$ (where 
$C(X)$ is the algebra of continuous complex 
functions on $X$). Letting $I$ be the closed ideal consisting of elements of 
$A$ which vanish on $X$, the least constant $k$ 
such that $$k\sup\{|f(x)|:x\in X\}\ge \Vert f+I\Vert\ \hbox{ for all }f\in A$$ 
is called the {\sl Helson constant} of $X$. We 
say that a Banach function algebra $A$ has the {\sl Helson property} (with 
constant $K$) if there is a constant $K$ such that 
whenever $U$ is a non-empty Gelfand open subset of $Max(A)$ there is an 
increasing net $(F_{\alpha})_{\alpha}$ of Helson sets 
of constant bounded by $K$ contained in $U$ such that $\bigcup_{\alpha} 
F_{\alpha}$ is Gelfand dense in $U$. 

A subset $E$ of an abelian group $G$ is said to be an {\sl independent set} 
[\the\Rud; 5.1.1] if it has the following property 
(following [\the\Rud] we use additive notation): for every choice of distinct 
points $x_1,\ldots ,x_k$ of $E$ and integers 
$n_1,\ldots ,n_k$, either 
$$n_1x_1=n_2x_2=\ldots =n_kx_k=0$$ 
or $$n_1x_1+n_2x_2+\ldots +n_kx_k\ne0.$$ The importance of this definition is 
that if $G$ is a locally compact abelian group 
and $F\subseteq\widehat G$ is a finite independent set then $F$ is a Helson 
set with Helson constant bounded by $2$ [\the\Rud; 
5.6.7]. 

A locally compact abelian group is an {\sl I-group} [\the\Rud; 2.5.5] if every 
neighbourhood of the identity contains an 
element of infinite order. A simple argument shows that this implies that 
every non-empty open subset contains an element of 
infinite order---indeed the set of elements of finite order is meagre. The 
next lemma is similar to [\the\Rud; 5.2.3]. 
\bigskip 
\noindent {\bf Lemma 3.3} {\sl Let $N$ be a locally compact, second countable 
abelian group and suppose either that $N$ is an 
I-group or that $N$ is an uncountable group in which every non-trivial element 
has order $q$ with $q$ prime. 
Let $U$ be a non-empty open subset of $N$. Then there exists an increasing net 
$(F_{\alpha})$ of finite independent sets 
contained in $U$, such that $\bigcup_{\alpha}F_{\alpha}$ is dense in $U$.} 
\bigskip 
\noindent {\bf Proof.} Suppose first that $N$ is an I-group. Let $U$ be a 
non-empty open subset of $N$ and let $(X_i)$ be a 
base of non-empty open sets for the topology on $U$. Choose $x_1\in X_1$ of 
infinite order. Suppose, for an inductive 
hypothesis, that we have chosen an independent set $\{x_1,\ldots ,x_k\}$ with 
each $x_i\in X_i$ and of infinite order. Let $S$ 
be the countable set $S=\{ m_1x_1+\ldots +m_kx_k: m_i\in {\bf Z}\}$. Let $s\in 
S$ and let $m\in {\bf Z}\setminus \{ 0\}$. Let 
$sol_m(s)=\{ x\in X_{k+1}: mx=s\}$. Then either $sol_m(s)$ is empty (hence 
meagre) or else there exists $y\in sol_m(s)$. In 
this case $$\eqalign{sol_m(s)&=\{ z\in X_{k+1}: m(y-z)=0\}\cr 
&=\{ x\in N: y+x\in X_{k+1},\ mx=0\}.\cr}$$ But $\{ x\in N:mx=0\}$ is a closed 
set, with no interior since $N$ is an I-group, 
so once again $sol_m(s)$ is meagre. Thus $\bigcup\{ sol_m(s):s\in S, m\in {\bf 
Z}\setminus\{0\}\}$ 
is also meagre, so $$\bigcup\{ sol_m(s):s\in S, m\in {\bf Z}\setminus\{0\}\} 
\cap X_{k+1}$$ is meagre in $X_{k+1}$. Since the set of elements of finite 
order is also meagre, there exists $x_{k+1}\in 
X_{k+1}\backslash 
\bigcup\{ sol_m(s):s\in S, m\in {\bf Z}\setminus\{0\}\}$ with $x_{k+1}$ of 
infinite order. It is straightforward to check that 
$\{ x_1,\dots, x_k,x_{k+1}\}$ satisfies the inductive hypothesis. Evidently 
$\bigcup _{i=1}^{\infty} \{x_i\}$ is dense in $U$. 

Now suppose that $N$ is uncountable and that every non-trivial element of $N$ 
has order $q$ with $q$ prime. As before, let $U$ 
be an open subset of $N$ and let $(X_i)_{i\ge 1}$ be a base of non-empty open 
sets for the topology on $U$. Choose $0\ne x_1\in 
X_1$. Suppose, for an inductive hypothesis, that we have chosen an independent 
set $\{x_1,\ldots ,x_k\}$ with each $0\ne x_i\in 
X_i$. Let $S$ be the countable set $S=\{ m_1x_1+\ldots +m_kx_k: m_i\in {\bf 
Z}\}$. Let $s\in S$ and let $m\in {\bf Z}\setminus 
q{\bf Z}$. Let $sol_m(s)=\{ x\in X_{k+1}: mx=s\}$. Then either $sol_m(s)$ is 
empty or else there exists $y\in sol_m(s)$. In 
this case $$sol_m(s)=\{ z\in X_{k+1}: m(y-z)=0\}=\{ y\}$$ since $q$ does not 
divide $m$. Thus $\bigcup\{ sol_m(s):s\in S, m\in 
{\bf Z}\setminus q{\bf Z}\}$ 
is finite or countably infinite, 
so $$X_{k+1}\setminus\bigcup\{ sol_m(s):s\in S, m\in {\bf Z}\setminus q{\bf 
Z}\}$$ is uncountable. Let 
$0\ne x_{k+1}\in X_{k+1}\backslash 
\bigcup\{ sol_m(s):s\in S, m\in {\bf Z}\setminus q{\bf Z}\}$. It is 
straightforward to check that $\{ x_1,\dots, x_k,x_{k+1}\}$ 
satisfies the inductive hypothesis. Evidently $\bigcup _{i=1}^{\infty} 
\{x_i\}$ is dense in $U$. Q.E.D. 
\bigskip 
\noindent It follows from Lemma 3.3 that if $N$ is a locally compact abelian 
group such that $\widehat N$ is either a second 
countable I-group or a second countable, uncountable group in which every 
non-trivial element has order $q$ with $q$ prime 
then $L^1(N)$ has the Helson property (with constant $2$). In the case of the 
second countable I-groups, the independent set 
$F$ constructed in Lemma 3.3 consists of elements of infinite order. This 
implies that $F$ is actually a Helson set with Helson 
constant $1$ [\the\Rud; 5.1.3 Corollary, and 5.5.2]. 

For the next lemma, let $D_q$ be the compact group obtained as the direct 
product of countably many copies of the cyclic group 
of order $q$, where $q$ is an integer, $q\ge 2$. 
\bigskip 
\noindent{\bf Lemma 3.4} {\sl Let $G$ be a locally compact group with an 
abelian, non-compact, closed, normal subgroup $N$ of 
finite index $m$. Then $G$ has a quotient $G'$ which has a non-compact closed 
normal subgroup $N'$ of finite index, such that 
$\widehat{N'}$ is either a second countable I-group or a second countable, 
uncountable group in which every non-trivial element 
has order $q$ where $q$ is a prime.} 
\bigskip 
\noindent{\bf Proof.} 
We begin by reducing to the case where $G$ is discrete.
Since $G$ has an abelian subgroup of finite index,
$G$ is a projective limit of Lie groups. Thus, passing to
a quotient modulo some compact normal subgroup, we may assume that
$G$ is a Lie group. Let $N_0$ denote the connected component of $N$
containing the identity. Then $N_0$ is open in $G$.

Suppose first that $N/N_0$ is finite. Then, by the structure theory of locally
compact, abelian groups, $N$ is isomorphic to ${\bf R}^n \times M$ where $M$ 
is a compact group and $n \geq 1$ (since $N$ is non-compact). Then $M$ is normal
in $G$, and $\widehat{N/M}$ is isomorphic to ${\bf R}^n$, which is a second countable 
I-group. Thus we are left with the case where $N/N_0$ is infinite, and so, 
passing to $G/N_0$ we may assume that $G$ is discrete.

By [\the\Rud; 2.5.5] there is a closed subgroup $H$ of 
$N$ such that either $\widehat {N/H}$ is a second 
countable I-group or such that $\widehat {N/H}$ is isomorphic to $D_q$ for 
some prime $q$. However $H$ may not be normal in 
$G$, so set $\widetilde H=\bigcap _{x\in G} x^{-1}Hx$. Then $\widetilde H$ is 
normal in $G$. Choose coset representatives 
$x_1,\ldots ,x_m$ for $G/N$, and for $1\le i\le m$ set $H_i=x_i^{-1}Hx_i$. 
Then the homomorphism 
$$\widehat{N/H_1}\times\cdots\times \widehat{N/ H_m}\to\widehat{N/\widetilde 
H}$$ from the product of the subgroups 
$\widehat{N/ H_i}$ of $\widehat{N/\widetilde H}$ given by 
$(\alpha_1,\ldots,\alpha_m)\mapsto \alpha_1\ldots\alpha_m$ is 
continuous and has dense range in $\widehat {N/\widetilde H}$, since this 
range separates the points in $N/{\widetilde H}$.  Since $N$ is discrete, all 
of the $\widehat{N/H_i}$ are compact and hence the above homomorphism is
surjective. Thus $\widehat{N/{\widetilde H}}$ is a quotient of the product
$\widehat{N/H_1}\times\cdots\times \widehat{N/ H_m}$.
Note also that $\widehat {N/\widetilde H}$ is compact and infinite, so is 
certainly uncountable.

Now if $\widehat{N /H}$ is isomorphic to $D_q$, then every 
$\widehat{N/H_i}$ is isomorphic to $D_q$ and this implies that 
$\widehat{N/{\widetilde H}}$ is second countable, and that
every non-trivial element of $\widehat{N/{\widetilde H}}$ has order $q$. 
Finally, if $\widehat {N/H}$ is a second countable I-group, then 
$\widehat{N/H_1}\times\cdots\times \widehat{N/ H_m}$ is also 
a second countable I-group. Clearly, then, $\widehat{N/{\widetilde H}}$ is
second countable. Since the homomorphism above is continuous and each 
$\widehat{N/H_i}$ is an I-group which is a subgroup of 
$\widehat{N/{\widetilde H}}$, the latter is an I-group as well.
Q.E.D. 
\bigskip 
\noindent We now introduce some further notation. For a locally compact group 
$G$ and a closed subset $X$ of $\widehat G$, let 
$K(X)$ denote the kernel of $X$ in $L^1(G)$. Then $K(X)$ is the largest ideal 
of $L^1(G)$ with hull equal to $X$. Ludwig has 
shown that if $G$ has polynomial growth and $L^1(G)$ is hermitian then there 
is a smallest closed ideal $j(X)$ whose hull is 
equal to $X$ [\the\Lud]. 

Now let $G$ be a locally compact group and suppose that $G$ has an abelian, 
closed, normal subgroup $N$ of finite index. Let 
$G$ act on $\widehat N$ in the usual way by $(x,\alpha)\mapsto x\cdot \alpha$ 
$(x\in G,\alpha\in\widehat N)$, where 
$x\cdot\alpha(n)=\alpha (x^{-1}nx)$ $(n\in N)$. 
Let $\widehat N /G$ be the set of $G$-orbits in $\widehat N$, with the 
quotient topology, and let $q:\widehat N\to \widehat N 
/G$ be the quotient map. Let $\phi:\widehat G\to\widehat N/G$ be the map 
defined $\phi(\pi)=G\cdot\alpha$, where 
$\alpha\in\widehat N$ and $\ker(\alpha)\supseteq\ker(\pi|_N)$. Then $q$ and 
$\phi$ are both continuous and both open. For a set 
$C\subseteq \widehat N$, let $\widetilde C=\phi^{-1}(q(W))$. 
\bigskip 
\noindent {\bf Lemma 3.5} {\sl Let $G$ be a locally compact group and suppose 
that $G$ has an abelian, closed, normal subgroup 
$N$ of finite index. Suppose that $X$ is a $G$-invariant compact subset of 
$\widehat N$ and that $Y$ is a $G$-invariant closed 
subset of $\widehat N$ such that $X$ is contained in the interior of $Y$ and 
such that the complement ot $Y$ in $\widehat N$ is 
relatively compact. Then $K(\widetilde Y)\subseteq j(\widetilde X)$.} 
\bigskip 
\noindent {\bf Proof.} Let $f\in K(\widetilde Y)$. By [\the\Lud] it is enough 
to find $g\in K(\widetilde X)$ such that $gf=f$. 
Let $h\in L^1(N)$ such that $h\in K(X)$ and $\hat h$ takes the constant value 
$1$ on the closure of the complement of $Y$. Such 
a function exists because $L^1(N)$ is completely regular and the complement of 
$Y$ is relatively compact. Since $N$ is an open 
subset of $G$ we may extend $h$ to an element $g\in L^1(G)$ by setting 
$h(x)=0$ for $x\in G\backslash N$. The element $g$ has 
the required property. Thus $f\in j(\widetilde X)$. Q.E.D. 

\bigskip 
\noindent Part of the argument of the next theorem adapts a method used in 
[\the\FS; 1.2]. 
\bigskip 
\noindent {\bf Theorem 3.6} {\sl Let $G$ be a locally compact group and 
suppose that $G$ has an abelian, non-compact, closed, 
normal subgroup $N$ of finite index $m$. Then $L^1(G)$ does not have spectral 
synthesis, and $\tau_r$ is not Hausdorff on 
$Id(L^1(G))$.} 
\bigskip 
\noindent {\bf Proof.} The properties of spectral synthesis and of $\tau_r$ 
being Hausdorff both pass to quotients [\the\Id; 
2.9]. Thus by Lemma 3.4 we may suppose that the dual $\widehat N$ of the 
abelian subgroup $N$ is either a second countable 
I-group or a second countable, uncountable group in which every non-trivial 
element has order $q$, with $q$ prime. (This 
reduction step is not required for the proof that $L^1(G)$ does not have 
spectral synthesis). 

For $\alpha\in\widehat N$, let $G_{\alpha}$ denote its stabilizer in $G$. The 
lengths of orbits of elements of $\widehat N$ are 
bounded by $m$. Let $d$ be the maximal orbit length, and set $\widehat 
N_d=\{\alpha\in\widehat N :|G:G_{\alpha}|=d\}$. Then 
$\widehat N_d$ is non-empty and open. Let ${\cal H}$ be the set of subgroups 
of $G$ of index $d$. Then $\widehat 
N_d=\bigcup_{H\in{\cal H}}\widehat N_H$, where the union is disjoint, and 
$\widehat N_H:=\{\alpha\in\widehat N: 
G_{\alpha}=H\}$. 

Fix $H\in{\cal H}$. Then there exists a non-empty open set $W\subseteq 
\widehat N_H$ such that the sets $x\cdot W$ $(x\in G/N)$ 
are pairwise disjoint. By Proposition 3.2, $W$ contains a non-spectral set $F$ 
for $L^1(N)$. Set $X=G\cdot F=\{x\cdot F:x\in 
G/N\}$. Then $X$ is $G$-invariant, and is also non-spectral since it has $F$ 
as a non-spectral non-empty clopen subset. Hence 
$\widetilde X$ is also non-spectral by the projection theorem [\the\HaLu; 
2.6], which we may use because $G$ has polynomial 
growth and $L^1(G)$ is hermitian. This shows that $L^1(G)$ does not have 
spectral synthesis. 

Now let $(V_{\alpha})_{\alpha}$ be a net of $G$-invariant, decreasing, open 
neighbourhoods of $X$, each having compact 
complement in $\widehat N$, such that $\bigcap_{\alpha} N_{\alpha}=X$ (where 
for each $\alpha$, $N_{\alpha}$ is the closure of 
$V_{\alpha}$). Then $(\widetilde V_{\alpha})_{\alpha}$ is a net of decreasing 
open subsets of $\widehat G$, and 
$\bigcap_{\alpha} \widetilde N_{\alpha} =\widetilde X$. Hence $(K(\widetilde 
N_{\alpha} ))_{\alpha}$ is an increasing net in 
$Id(L^1(G))$, and $K(\widetilde N_{\alpha})\subseteq j(\widetilde X)$ for all 
$\alpha$ by Lemma 3.5. Hence 
$$I:=\overline{\bigcup_{\alpha} K(\widetilde N_{\alpha})}\ \subseteq 
j(\widetilde X).$$ 

By Lemma 3.3 there is for each $\alpha$ an increasing net 
$(H_{\beta(\alpha)})_{\beta(\alpha)}$ of Helson sets in $V_{\alpha}$, 
each of Helson constant bounded by $2$, such that $\bigcup_{\beta(\alpha)} 
H_{\beta(\alpha)}$ is dense in $V_{\alpha}$. Then 
$F_{\beta(\alpha)}:=G\cdot H_{\beta(\alpha)}$ is a $G$-invariant Helson set, 
with a bound $K$ dependent on $m$, but independent 
of $\alpha$, and of course $\bigcup_{\beta(\alpha)} F_{\beta(\alpha)}$ is also 
dense in $V_{\alpha}$. Hence 
$(K({\widetilde {F_{\beta(\alpha)}}} ))_{\beta(\alpha)}$ is a decreasing net 
in 
$Id(L^1(G))$, and $\bigcap_{\beta(\alpha)}K({\widetilde {F_{\beta(\alpha)}}} ) 
= 
K(\widetilde N_{\alpha} )$. Hence $K({\widetilde {F_{\beta(\alpha)}}} 
)\mapright{\beta(\alpha)} 
K(\widetilde N_{\alpha} )$ $(\tau_r)$ by Lemma 1.1. But $K(\widetilde 
N_{\alpha} )\to 
I$ $(\tau_r)$, also by Lemma 1.1, so if $(K(\widetilde {F_{\gamma}} 
))_{\gamma}$ denotes the `diagonal' net, see [\the\Kel; \S 
2, Theorem 4], then $K(\widetilde {F_{\gamma}} )\to I$ $(\tau_r)$. 

Now we show that $K(\widetilde F_{\gamma} )\to K(\widetilde X)$ $(\tau_r)$. 
First we establish convergence for $\tau_n$. 
Suppose that $f\notin K(\widetilde X)$. Then there exists $P\in \widetilde X$ 
such that $f\notin P$. For each $Q\in 
Prim_*(L^1(G))$, let $Q_*$ be the primitive ideal of $C^*(G)$ whose 
intersection with $L^1(G)$ is equal to $Q$. Set 
$\epsilon=\frac12 \Vert f+P_*\Vert_*$, where $\Vert\,\cdot\,\Vert_*$ denotes 
the C$^*$-norm on $C^*(G)$. By the lower 
semicontinuity of norm functions for C$^*$-algebras, there is a hull-kernel 
neighbourhood $M_*$ of $P_*$ in $Prim(C^*(G))$ such 
that $\Vert f+Q_*\Vert_*>\epsilon$ for all $Q_*\in M_*$. Let $M=\{ Q\in 
Prim_*(L^1(G)):Q_*\in M_*\}$. Then $M$ is an open 
neighbourhood of $P$ in $Prim_*(L^1(G))$, because $G$ is $^*$-regular, so 
eventually there is for each $\gamma$ and element 
$Q_{\gamma}\in\widetilde F_{\gamma}\cap M$. Hence eventually 
$$\Vert f+K(\widetilde F_{\gamma})\Vert\ge\Vert f+Q_{\gamma}\Vert\ge\Vert 
f+(Q_{\gamma})_*\Vert_*>\epsilon.$$ This shows that 
$K(\widetilde F_{\gamma} )\to K(\widetilde X)$ $(\tau_n)$. 

Finally we show that $K(\widetilde F_{\gamma} )\to K(\widetilde X)$ 
$(\tau_u)$. By [\the\Id; 2.1] it is enough to show that for 
each $f\in K(\widetilde X)$ and $\epsilon>0$ there is a neighbourhood 
$\widetilde M$ of $\widetilde X$ in $\widehat G$ such 
that $\Vert f+K(\widetilde F_{\gamma})\Vert<\epsilon$ whenever $\widetilde 
F_{\gamma}\subseteq \widetilde M$. Since $N$ is an 
open subgroup of $G$ there is a projection $P:L^1(G)\to L^1(N)$ given by 
$P(g)=g|_N$ $(g\in L^1(G))$. Then for each coset 
representative $x\in G/N$, $P(L_xf)\in L^1(N)$, and in fact $P(L_xf)\in K(F)$. 
Thus there exists a neighbourhood $M_x$ of $F$ 
in $\widehat N$ such that $|P(L_xf)(s)|<\epsilon/(Km)$ for each $s\in M_x$. 
Hence by the Helson property, $\Vert P(L_xf) 
+K(F_{\gamma})\Vert<\epsilon/m$ whenever $F_{\gamma}\subseteq M_x$. Set 
$M=\bigcap_{x\in G/N} M_x$. Then whenever 
$F_{\gamma}\subseteq M$, there is, for each coset representative $x\in G/N$, 
an $h_x\in K(F_\gamma)$ such that $\Vert 
P(L_xf)-h_x\Vert <\epsilon/m$. 
Define $h\in K(\widetilde F_{\gamma})$ by $h(xn)=h_x(n)$ $(n\in N)$. Then 
$$\eqalign {\Vert f-h\Vert=&\left\Vert\sum_{x\in G/N} f 1_{xN} 
-h1_{xN}\right\Vert\cr 
=&\sum_{x\in G/N} \Vert f1_{xN}-h1_{xN}\Vert\cr 
=&\sum_{x\in G/N} \Vert P(L_xf)-h_x\Vert<\epsilon.\cr 
}$$ 
Hence 
$K(\widetilde F_{\gamma})\to K(\widetilde X)$ $(\tau_u)$, so 
$$K(\widetilde F_{\gamma})\to K(\widetilde X)~~(\tau_r).$$ 
Since $I\subseteq j(\widetilde X)$ 
and $j(\widetilde X)$ is a strict subset of 
$K(\widetilde X)$, we have $I\ne K(\widetilde X)$, 
so $\tau_r$ is not Hausdorff. Q.E.D. 
\bigskip 

\noindent Note, for the next result, that if a group $G$ has a compact normal 
subgroup $K$ such that $G/K$ is a finite 
extension of a nilpotent group then $L^1(G)$ is hermitian and $G$ has 
polynomial growth [\the\Palm], and hence $L^1(G)$ is 
$^*$-regular [\the\Bar]. Thus Theorem 2.1 applies. 
\medskip 
\noindent {\bf Theorem 3.7} {\sl Let $G$ be a locally compact group. Suppose 
that $G$ has a compact normal subgroup $K$ such 
that $G/K$ is a finite extension of a nilpotent group. Then the following are 
equivalent: 

(i) $\tau_r$ is Hausdorff on $Id(L^1(G))$, 

(ii) $\tau_{\infty}$ is Hausdorff on $Id(L^1(G))$, 

(iii) $\tau_r$ and $\tau_{\infty}$ coincide on $Id(L^1(G))$, 

(iv) $G$ is compact, 

(v) $L^1(G)$ has spectral synthesis.} 
\bigskip 
\noindent {\bf Proof.} (iv)$\Rightarrow$(ii) is proved in [\the\Bec; p.72]. 
The equivalence of (ii) and (iii), and hence the 
implication (ii)$\Rightarrow$(i), is proved in [\the\Id; 3.1.1]. The 
implication (iv)$\Rightarrow$(v) is well-known. It remains 
to prove that (i)$\Rightarrow$(iv) and that (v)$\Rightarrow$(iv). First we 
show that (i)$\Rightarrow$(iv). 

Suppose that $\tau_r$ is Hausdorff on $Id(L^1(G))$, and suppose too, to begin 
with, that $G$ is a finite extension of an 
[FD]$^-$ group, so that $G$ has a closed normal subgroup $N$ of finite index 
with $N\in$[FD]$^-$. Let $K$ be the closure of the 
commutator subgroup of $N$. Then $K$ is compact, and $K$ is normal in $G$. 
Hence $G/K$ is a finite extension of the abelian 
group $N/K$, and $\tau_r$ is Hausdorff on $Id(L^1(G/K))$. Thus $G/K$ is 
compact by Theorem 3.6, so $G$ is compact. 

Now we consider the general case when $G$ has a compact normal subgroup $K$ 
such that $G/K$ is a finite extension of a 
nilpotent group. It is enough to show that $G/K$ is compact, and hence we may 
suppose that $G$ itself has a nilpotent closed 
normal subgroup $N$ of finite index. We proceed by induction on the length 
$l(N)$ of nilpotency of $N$. If $l(N)=1$ then $N$ is 
abelian, so it follows from Theorem 3.6 that $G$ is compact. Now suppose that 
we have established the result for $l(N)=k\ge 1$, 
and that $l(N)=k+1$. Let $Z(N)$ be the centre of $N$. Then $Z(N)$ is normal in 
$G$, so $G/Z(N)$ is a finite extension of 
$N/Z(N)$, and $l(N/Z(N))=k$, so $G/Z(N)$ is compact. Hence $N/Z(N)$ is 
compact, so $N\in$[Z], i.e. $N$ is a central group. But 
[Z]$\subseteq$[FD]$^-$ [\the\GM], so we have that $G$ is a finite extension of 
an [FD]$^-$-group. Hence $G$ is compact by the 
previous paragraph. 

Finally we show that (v)$\Rightarrow$(iv). This follows by exactly the same 
argument used for (i)$\Rightarrow$(iv), replacing 
the hypothesis that $\tau_r$ is Hausdorff by the hypothesis that $L^1(G)$ has 
spectral synthesis. This property also passes to 
quotients, and the appeals to Theorem 3.6 are still valid. Q.E.D.

\bigskip 
\noindent Theorem 3.7, of course, covers the classes of nilpotent groups and 
[FD]$^-$-groups. Every Moore group is a finite 
extension of an [FD]$^-$-group [\the\Robe], so Theorem 3.7 also covers Moore 
groups. 

\bigskip 
\noindent For the final result of the paper, we apply Theorem 3.5, in the 
abelian case, to the situation of Fourier algebras. 
\bigskip 
\noindent {\bf Theorem 3.8} {\sl Let $G$ be a non-discrete locally compact 
group with Fourier algebra $A(G)$. Then $\tau_r$ is 
not Hausdorff on $Id(A(G))$.} 
\bigskip 
\noindent {\bf Proof.} Suppose that $\tau_r$ is Hausdorff on $Id(A(G))$. Then 
$\tau_r$ is Hausdorff on $Id(A(H))$ for every 
closed subgroup $H$ of $G$, because $A(H)$ is isomorphic to $A(G)/I(H)$, where 
$I(H)$ is the ideal consisting of all functions 
in $A(G)$ which vanish on $H$ [\the\For; Lemma 3.8]. 

Suppose first that $G$ is totally disconnected. Then because $G$ is not 
discrete, $G$ has an infinite, compact open subgroup 
$K$. By [\the\Zelm; Theorem 2], $K$ has an infinite closed abelian subgroup 
$H$. Then $\tau_r$ is Hausdorff on $Id(A(H))$, and 
$A(H)=L^1(\widehat H )$, so $\widehat H$ is compact by Theorem 3.6. Hence $H$ 
is discrete, contradicting the fact that it is 
infinite and compact. 

It remains therefore to show that if $\tau_r$ is Hausdorff on $Id(A(G))$ then 
$G$ must be totally disconnected. Suppose, for a 
contradiction, that $G_0$, the connected component of the identity, is 
non-trivial. A connected group is a projective limit of 
Lie groups, so there is a compact normal subgroup $K$ of $G_0$ such that 
$G_0/K$ is a non-trivial Lie group. But $A(G_0/K)$ is 
a quotient of $A(G_0)$. Indeed the mapping $u\mapsto \dot u$, where $$\dot u 
(xK)=\int_K u(xt) dt$$ (with normalized Haar 
measure on $K$) is a continuous homomorphism of $A(G_0)$ onto $A(G_0/K)$ 
[\the\Ey]. Thus $\tau_r$ is Hausdorff on 
$Id(A(G_0/K))$. On the other hand, $G_0/K$ is a connected Lie group, so it is 
generated by its one parameter subgroups, i.e. 
images of analytic homomorphisms of ${\bf R}$ into $G_0/K$. Hence $G_0/K$ 
contains numerous non-discrete closed abelian 
subgroups. As in the previous paragraph, this contradicts the fact that 
$\tau_r$ is Hausdorff on $Id(A(F))$ for each closed 
abelian subgroup $F$ of $G_0/K$, which forces each such $F$ to be discrete. 
Q.E.D. 
\bigskip 
\noindent It was shown in [\the\KL] that if $G$ is a discrete group then 
$A(G)$ has spectral synthesis (and hence 
$\tau_{\infty}$ and $\tau_r$ are Hausdorff) provided that $A(G)$ satisfies an 
additional weak condition. This additional 
condition is satisfied whenever $G$ is an amenable discrete group, and is 
probably satisfied for all discrete groups. 

\bigskip 
We close this section by observing that the Hausdorffness of $\tau_r$ is {\sl 
not} equivalent to spectral synthesis for general 
Banach $^*$-algebras. The Banach algebra $C^1[0,1]$ is a hermitian, (locally) 
regular, $^*$-semisimple Banach $^*$-algebra 
without spectral synthesis, but $\tau_r$ is Hausdorff on $Id(C^1[0,1])$ 
[\the\FS]. 

\bigskip 
\bigskip 
\centerline{\bf References} 
\medskip 
\input TempReferences.tex 
\bigskip 
\bigskip 
\centerline{School of Mathematical Sciences} 
\centerline{University of Nottingham} 
\centerline{ NG7 2RD} 
\centerline{ U.K.} 
\medskip 
\centerline{email: Joel.Feinstein@nottingham.ac.uk} 
\bigskip 
\centerline{Fachbereich Mathematik/Informatik} 
\centerline{Universit\"at Paderborn} 
\centerline{D-33095 Paderborn} 
\centerline{Germany} 
\medskip 
\centerline{email: kaniuth@uni-paderborn.de} 
\bigskip 
\centerline{Department of Mathematical Sciences} 
\centerline{University of Aberdeen} 
\centerline{AB24 3UE} 
\centerline{U.K.} 
\medskip 
\centerline{e-mail: ds@maths.abdn.ac.uk} 
\end